\documentclass [a4paper 12pt] {article}

\usepackage[a4paper]{geometry}
\usepackage{pdflscape}
\usepackage[english]{babel}
\usepackage{amsmath}
\usepackage[utf8x]{inputenc}
\usepackage{siunitx}
\usepackage{setspace}
\usepackage{color}
\usepackage{graphicx}
\usepackage{hyperref}
\usepackage{cleveref}
\usepackage{enumitem}
\usepackage{multirow}
\usepackage[table]{xcolor}
\usepackage{chemfig}
\usepackage[standard]{ntheorem}

\begin{document}
\begin{center}

\begin{title}
\textbf{\textbf{A simplified expression for the solution of cubic polynomial equations using function evaluation}}\par
\vspace{0.3cm}
\textit {Ababu T. Tiruneh, University of Eswatini, Department of Env. Health Science, Swaziland.}\par
Email: \href{mailto:ababute@gmail.com}{\nolinkurl ababute@gmail.com}\par
\end{title}

\end{center}

\vspace{18pt}

\renewcommand{\abstractname}{Abstract \hfill}

\begin{abstract}

This paper presents a simplified method of expressing the solution to cubic equations in terms of function evaluation only. The method eliminates the need to manipulate the original coefficients of the cubic polynomial and makes the solution free from such coefficients. In addition, the usual substitution needed to reduce the cubics is implicit in that the final solution is expressed in terms only of the function and derivative values of the given cubic polynomial at a single point. The proposed methodology simplifies the solution to cubic equations making them easy to remember and solve.\\

 
\end{abstract}

\section{Introduction}

Polynomials of higher degree arise often in problems in science and engineering. According to the fundamental theorem of algebra, a polynomial equation of degree n has at most n distinct solutions \cite {Fathi}. The history of symbolic manipulation for solving polynomial equation notes the work of Luca Pacioli in 1494  \cite{Dorsey} in which the basis was laid for solving linear and quadratic equations while he stated the cubic equation as impossible to solve and asking the Italian mathematical community to take the challenge. Mathematicians seemed to have felt hitting a wall after solving quadratic equations as it took quite a while for the solution to cubic equation to be found \cite{Robert}. The solution to the cubic in the depressed form $x^3 + bx +c $ was discovered by del Ferro initially but he passed it to his student instead of publishing it. Eventually Caradano in 1540 got hold of the solution and published the results, crediting Del Ferro for the solution. Cardano also published Ferrari's solution of the quartic equation. \par

\vspace{12pt}

The general cubic of the form  $$ x^3 +bx^2 +cx +d $$ can be reduced to a depressed form by a suitable substitution involving a new variable\cite{Dickson1}. Accordingly,  the substitution takes the form of: $$ x = y - \frac{b}{3} $$

The depressed cubic takes the form: $$ y^3 + py + q$$  Where  p and q are expressed in the forms: $$ p = c- \frac{b^2}{3}  \quad   and  \quad q = d - \frac{bc}{3} + \frac{2b^3}{27}$$

Further substitution of the form: $$ y = z - \frac{p}{3z}$$ will reduce the depressed cubic after simplification to the form $$ z^6 + qz^3 - \frac{p^3}{27} = 0 $$  Solving this as a quadratic for $z^3$ gives: $$ z^3 = -\frac{q}{2} \pm \sqrt{R}  \quad  R= {\left(\frac{p}{3}\right)}^3 + {\left(\frac{q}{2}\right)}^2$$  \par

Considering the three cube roots of unity \cite{Dickson2} namely; 
$$ 1, \quad  \omega= -\frac{1}{2} + \frac{1}{2} \sqrt{3}i \quad  and \quad  {\omega}^2 = -\frac{1}{2} - \frac{1}{2}\sqrt{3}i$$

It can then be easily verified \cite{Dickson2} that the three solutions of the depressed cubic equations are given by:

$$ y_1 = A+ B , \quad   y_2 = \omega A + {\omega}^2 B   \quad   and \quad y_3 =\omega B + {\omega}^2 A$$

The value of A and B in the above equation are given by:

$$ A = \sqrt[3]{ -\frac{q}{2} + \sqrt{R}}  \quad   and \quad  B = \sqrt[3]{ -\frac{q}{2} - \sqrt{R}} $$

\vspace{12pt}

Other similar methods also arrived essentially at the same formula state above. For example see the work of  Mukundan \cite{Mukundan} that transforms the original depressed cubic equation into two cubic equivalents that are solved by simply taking the cube roots. A more comprehensive coverage of the formulae for the solution of the cubic equation is found at the Wolfram Mathworld website \cite{Wolfram} in which some of the formula used in the examples used in this paper are referred to. Lagrange in 1770-1771, as part of his study of higher degree equation used Fourier transform and subsequent inversion \cite{Lagrange}. He hoped to extend his methods to arbitrarily high degree polynomials. However, Lagrange noted that the resulting resolvent polynomial for a five degree polynomial equation was actually a six degree polynomial, prompting a possible hint to him that such equations may not be solvable. \par

\vspace{12pt} 

The  solution to cubic equation has laid a basis for the methods developed for solving quartic equations. For example, Leonard Euler (1707-1783) appreciated the central role of the resolvent cubic in the solution of the quartic polynomial equations \cite{Nickalls}. Euler \cite{Euler}, through his solution to the depressed (in which the $x^3$ term is zero) quartic equations, showed that each root is expressed as the sum of three square roots that are solutions of the resolvent cubic equation. Even the first solution of the quartic equation by Del Ferro that was published together with the solution of the cubic equation by Cardano involved converting the original quartic equation in to two complete squares. During the process of solving such equation a resolvent cubic equations arises that must be solved first \cite{Franz}.  Several of the recent methods for solving quartic equations also do require solving the resolvent cubic. For example see the works of  Saghe \cite{Saghe} Fathi \cite{Fathi2} and that of Kulkarni \cite{kulkarni}. \par 

\vspace{12pt}

 Solution to five degree polynomials and greater were met with limited success as it later transpired to mathematicians that the solutions could not be expressed in terms of radicals. The first clear proof came from Abel in 1824 who proved that the general polynomial of degree five can not be solved in terms of radicals. Abel later revised his proof by verifying that certain forms of the quintic equation (that were later termed as 'Abelian Galois group') can be solved by radicals. Galois in 1832 gave a proof that a five degree polynomial equation can be solved by radicals if and only if its Galois group is solvable \cite{Franz}.Certain polynomials equations of higher degree also  require solving a resolvent cubic. Example is the so called palindromic polynomials with symmetric coefficients. A six degree polynomial  that is palindromic results in a resolvent cubic through substitution of the fom $   y = x +\frac{1}{x} $  \cite{Franz}  \par

\vspace{12pt}

Some authors noted the relation of the transformation variables to the depressed cubic and quartic equations with derivative values. For example, Das, \cite{Das},  noted that the transformation to the depressed form of the quadratic, cubic and quartic equations correspond to the equations that set to zero the first derivative, second and third derivative respectively of the original equations. However, beyond this the transformation of the other coefficients of the cubic equation to functions and derivative values is not explained.  \par

\vspace{12pt}

\section{Methods}

Consider the cubic polynomial equation of the  form $ax^3 + bx^2 + cx + d =0$ Without loss of generality assume a=1 in which the equation reduces to  $x^3 + bx^2 + cx + d =0$  The case in which $a \neq 1$ will be dealt with a  simple revision of the solution at a later stage.\par
\vspace{12pt}

The methodology of arriving at the final simplified expression in terms of function evaluation follows the well known procedure of reducing the cubic to the deflated form and using Viete's magic substitution to reduce the equation to quadratic form. It should be known that in the end both procedures become implicit in the solution whereby the final expression consists of function value of the original cubic and its derivatives. \par
\vspace{12pt}

Define new variable t and constant z  such that \begin {equation} \label{eq.1} x = z + t \end{equation}  Substituting this new expression  from eq.1 into the original cubic equation $x^3 + bx^2 + cx + d =0$  results in the following:\par

$$x^3 + bx^2 + cx + d = {(z+t)}^3 + b{(z+t)}^2 + c(z+t) + d = 0$$ \par

After suitable simplification the above expression takes the form: 

\begin {equation} \label{eq.2}  t^3 +  (b+3z)t^2 + (3z^2 + 2bz + c)t+(z^3 +bz^2 + cz+d) = 0\end{equation}

\vspace{12pt}

It is clear that eq.2 contains expressions in the brackets that are functional and derivative values of the constant z defined in eq.1. Therefore, eq.2 can be expressed as:

\vspace{12pt}

\begin {equation} \label{eq.3}  t^3 +  \left(\frac{f''(z)}{2}\right)t^2 + (f'(z))t+f(z) = 0\end{equation}.

\vspace{12pt}

The deflated form is easily obtained by setting $ f''(z) = 6z+2b = 0 $   in eq.3 which gives $ z= -b/3$.  After this eq.3  will reduce to the reduced form given below.

\begin {equation} \label{eq.4}  t^3  + (f'(z))t + f(z) = 0\end{equation}

Now the famous Viete's substitution can be made in eq.4  by defining a variable s and a constant $\alpha$ such that: \begin{equation} \label{eq.5} t = s + \frac{\alpha}{s}\end {equation}

\vspace{12pt}

Substituting the expression in eq.5 into eq.4 and further simplification results in the following equation:

\vspace{12pt}

\begin {equation} \label{eq.6}  s^3  + \left(3\alpha +f'(z) \right)s + \left(3{\alpha}^2 + f'(z) \alpha \right) \frac{1}{s}  + f(z)  + \frac{{\alpha}^3}{s^3}= 0\end{equation}

\vspace{12pt}

It is apparent that the expression containing s and 1/s  in eq.6  will vanish if Viete's substitution, i.e., $ \alpha = -\frac{f'(z)}{3} $  is made.  After this substitution, eq.6 reduces to:

\vspace{12pt}

\begin {equation} \label{eq.7}  s^6  + f(z) s^3 - \frac{{f'(z)}^3}{27}= 0\end{equation}

\vspace{12pt}

Substituting $ r = s^3$  in eq.7  will reduce it to the quadratic form: \begin {equation} \label{eq.8}  r^2  + f(z) r - \frac{{f'(z)}^3}{27}= 0\end{equation}

\vspace{12pt}

The solution to eq.8 will be $$ r = \dfrac{-f(z) \pm \sqrt{{f(z)}^2 + 4\dfrac{{f'(z)}^3}{27}}}{2}$$

Substituting back $ s = r^{1/3}$ in the above expression results in the solution in terms of s, i.e., 

\vspace{12pt}

\begin {equation} \label{eq.9}  s = {\left[\dfrac{-f(z) \pm \sqrt{{f(z)}^2 + 4\dfrac{{f'(z)}^3}{27}}}{2}\right]}^{1/3}\end{equation}

\vspace{12pt}

Finally, substituting the expression in eq.9 in eq.5 and using the value of the constant z to reduce the cubic to the defalted form in eq.1 will give the solution of the cubic polynomial equation which will be in the form given below after simplification:

\vspace{12pt}

\begin {equation} \label{eq.10} x = {f''}^{-1}(0) +  {\left[\dfrac{-f(z) + \sqrt{{f(z)}^2 + 4\dfrac{{f'(z)}^3}{27}}}{2}\right]}^{1/3} + {\left[\dfrac{-f(z) - \sqrt{{f(z)}^2 + 4\dfrac{{f'(z)}^3}{27}}}{2}\right]}^{1/3}\end{equation}

\vspace{12pt}

It is shown, therefore, that the solution of the general cubic equation $x^3 + bx^2 + cx + d =0$ is expressed in  eq.1  in terms only of the the functional and derivative values at the constant z. The z value is found from the equation $f''(z) = 0$ and this value is given by $  z= {f''}^{-1}(0)$ in eq.10. It is apparent that eq.10 is free of the original coefficient of the cubic equation which might make it difficult to remember the solution in terms of these coefficients. The case of the general polynomial $ax^3 + bx^2 + cx + d =0$ is handled by dividing both  the functional value f(z) and the derivative f'(z) in eq.10 by the coefficient a which will eventually give the solution.     \par

Equation (10) can be reduced further with the following usual substitution of  cubic equations: $$   R= \frac{-f(z)}{2} \qquad and \qquad     Q= \frac{f'(z)}{3}$$

With the above definition of R and Q,  Equation (10) will be rewritten as:

\begin {equation} \label{eq.11} x = {f''}^{-1}(0) +  {\left[R+ \sqrt{{R}^2 + {Q}^3}\right]}^{1/3} + {\left[R- \sqrt{{R}^2 + {Q}^3}\right]}^{1/3}\end{equation}

The application of this simplified expression will be illustrated with the three examples given in the following section. \par

\vspace{12pt}

\begin{example}
 $ f(x) = x^3 - 6x^2 + 11x - 6 =0$  \par
\vspace{12pt}

It will be shown below that the all the roots of the above equations are real numbers. To begin with, we evaluate z, f(z)  and f'(z) \par

$f''(z) = 6z-12 =0$ means   $z= \frac{12}{6} = 2$
$$f'(z) = 3{z}^2-12(z)+11 = 3{(2)}^2-12(2)+11 = -1$$
$$f(z) ={2}^3-6({2}^2)+11(2)-6 = 0$$

$$   Q= \frac{f'(z)}{3} = \frac{-1}{3} \quad and \quad   R= \frac{-f(z)}{2} = 0$$

\vspace{12pt}

The discriminant D is evaluated as: $ D = {Q}^3 + {R}^2 = {(-\frac{1}{3})}^3 + {(0)}^2 = -\frac{1}{27} <0$   implying that all the roots are real and unequal since D \textless 0. The three roots are then computed as follows: \par

\vspace{12pt}

$$ \theta = {cos}^{-1}\left(\frac{R}{\sqrt{{(-Q)}^{3}}} \right) = {cos}^{-1} \left(\frac{0}{\sqrt{{(\frac{1}{3})}^{3}}} \right) =  {cos}^{-1} (0) = \frac{\pi}{2}$$

$$ x_1 = 2\sqrt{(-Q)} \quad cos (\frac{\theta}{3})+ z  = 2\sqrt{(\frac{1}{3})}\quad cos (\frac{\pi}{6})+ 2 = 1+2  = 3$$

$$ x_2 = 2\sqrt{(-Q)} \quad cos (\frac{\theta+ 2\pi}{3})+ z  = 2\sqrt{(\frac{1}{3})}\quad cos (\frac{5\pi}{6})+ 2 = -1+2  = 1$$

$$ x_3 = 2\sqrt{(-Q)} \quad cos (\frac{\theta+ 4\pi}{3})+ z  = 2\sqrt{(\frac{1}{3})}\quad cos (\frac{3\pi}{2})+ 2 = 0+2  = 2$$

\vspace{12pt}

Therefore, the roots of the given cubic equation are: $ x = \{1, 2, 3\}$  \par

\end{example}

\vspace{12pt}

\begin{example}
 $ f(x) = x^3 - 15x - 4 =0$  \par
\vspace{12pt}

This equation is used here for historical reasons as it was given by Bombeli in 1572. Again we begin by evaluating z, f(z)  and f'(z) \par

$f''(z) = 6z =0$ means   $z = 0$
$$f'(z) = 3{z}^2-15 = 3{(0)}^2-15 = -15$$
$$f(z) ={(0)}^3-15(0)-4 = -4$$

$$   Q= \frac{f'(z)}{3} = \frac{-15}{3} = -5 \quad and \quad   R= \frac{-f(z)}{2} = -\left(\frac{-4}{2}\right) =2 $$

\vspace{12pt}

The discriminant D is evaluated as: $ D = {Q}^3 + {R}^2 = {(-5)}^3 + {(2)}^2 = --125+4 = -121 <0$   implying again that all the roots are real and unequal since D \textless 0. The three roots are then computed as follows: \par

\vspace{12pt}

$$ \theta = {cos}^{-1}\left(\frac{R}{\sqrt{{(-Q)}^{3}}} \right) = {cos}^{-1} \left(\frac{2}{\sqrt{{(125)}}} \right) = 1.390942827$$

$$ x_1 = 2\sqrt{(-Q)} \quad cos (\frac{\theta}{3})+ z  = 2\sqrt{(5)}\quad cos (\frac{\theta}{3})+ 0= 4+0  = 4$$

$$ x_2 = 2\sqrt{(-Q)} \quad cos (\frac{\theta+ 2\pi}{3})+ z  = -2- \sqrt{3}$$

$$ x_3 = 2\sqrt{(-Q)} \quad cos (\frac{\theta+ 4\pi}{3})+ z  = -2+ \sqrt{3}$$

\vspace{12pt}

Therefore, the roots of the given cubic equation are: $ x = \{4, -2- \sqrt{3}, -2+\sqrt{3}\}$  \par

\end{example}

\vspace{12pt}

\begin{example}
 $ f(x) = x^3 - 5x^2 + 9x - 9 =0$  \par
\vspace{12pt}

It will be shown below that the roots of this equation include complex numbers. To begin with, we evaluate z, f(z)  and f'(z) \par
\vspace{12pt}

$f''(z) = 6z-10 =0$ means   $z= \frac{10}{6} = \frac{5}{3}$
$$f'(z) = 3{z}^2-10(z)+9 = 3{\left(\frac{5}{3}\right)}^2-10\left(\frac{5}{3}\right)+9 = \frac{2}{3}$$
$$f(z) ={\left(\frac{5}{3}\right)}^3-5{\left(\frac{5}{3}\right)}^2+9 \left(\frac{5}{3}\right)-9 = -\frac{88}{27}$$

$$   Q= \frac{f'(z)}{3} = \frac{\frac{2}{3}}{3} = \frac{2}{9} \quad and \quad   R= \frac{-f(z)}{2} = -\frac{1}{2} \left(\frac{-88}{27}\right) = \frac{44}{27}$$

\vspace{12pt}

The discriminant D is evaluated as: $ D = {Q}^3 + {R}^2 = {\left(\frac{2}{9}\right)}^3 + {\left(\frac{44}{27}\right)}^2  >0$   implying that the cubic equation has one real root and two complex conjugate roots.  The three roots are then computed as follows: \par

\vspace{12pt}

The real root  is found  from the equation: $ x_1 = z +B $   where B is given by:\par

\begin {equation} \label{eq.12} B =  {\left[R+ \sqrt{{R}^2 + {Q}^3}\right]}^{1/3} + {\left[R- \sqrt{{R}^2 + {Q}^3}\right]}^{1/3}\end{equation}

\vspace{12pt}

Substituting the computed values of R and Q in equation (12) gives $ B=  4/3$ \par

\vspace{12pt}

The real root is, therefore, $ x_1 =  z+B  = \frac{5}{3} + \frac{4}{3} = \frac{9}{3} = 3 $\par

\vspace{12pt}

The complex conjugate roots are obtained from the following equation:

\begin {equation} \label{eq.13} x_{2,3} =  z - \left(\frac{1}{2}\right)B \pm \left(\frac{\sqrt{3}}{2}\right)i \sqrt{{B}^2+4Q} \end{equation}

\vspace{12pt}

Substituting the computed values of $ z = \frac{5}{3} \quad  B = \frac{4}{3}  \quad  and \quad Q = \frac{44} {27} $  in equation (13) gives:\par

$$ X_{2,3} =  1 \pm \sqrt{2}i $$

\vspace{12pt}

The three roots of the given cubic equation are, therefore, $ x = \{3, \quad 1 + \sqrt{2}i, \quad 1 - \sqrt{2}i\}$  \par

\end{example}

\vspace{12pt}

\section{Conclusion}

The formula for the solution to cubic equations are a bit complex to work with. However, this paper demonstrated that it is possible to simplify the expression for the solutions of the cubic equations further through a method that requires only evaluation of the functional value and the derivative  at x value which makes the second derivative of the cubic equation zero.  When the solution is expressed essentially in terms of these  two  functional values, it would be easy to remember the formula for the solution as given by equation (10). In other words, the original general cubic equation containing the coefficients is solved in terms of functional values evaluated at a single point. There is no need to remember a formula containing the coefficients of the original cubic equation which can be complex for a general cubic. From educational perspective, this paper also provides an alternative mechanism for understanding and solving cubic equations using function evaluation. Such mechanism can also be extended to higher degree polynomial equations such as four degree polynomial and solvable quintics. The solution to quadratic equation using functional evaluation is, of course,  trivially simple and  likewise applicable.

\medskip
 
\bibliographystyle{unsrt}
\bibliography{samplebib_2}

\end{document}